\documentclass[12pt]{amsart}
\usepackage{amssymb,amscd}
\usepackage{verbatim}
\usepackage{graphicx}

\newtheorem*{Whitney towers}{Theorem~\ref{Whitney towers}}
\newtheorem*{h-towers}{Theorems ~\ref{half} \& \ref{$(n)$-solvable}}

\newtheorem*{surgery curves}{Theorem~\ref{surgery curves}}
\newtheorem*{cg=0}{Theorem~\ref{vanish}}

\newtheorem{thm}{Theorem}[section]

\newtheorem{prop}[thm]{Proposition}
\newtheorem{cla}[thm]{Claim}
\theoremstyle{definition}
\newtheorem{defn}[thm]{Definition}

\newtheorem{prob}[thm]{Problem}

\numberwithin{equation}{section}
\numberwithin{figure}{section}

\newcommand{\x}{\times}
\newcommand{\np}{\newpage}
\newcommand{\Z}{\mathbb{Z}}
\newcommand{\N}{\mathbb{N}}

\newcommand{\Q}{\mathbb{Q}}

\newcommand{\nl}{\newline}

\title{
Ribbon-moves of 2-knots: the torsion linking pairing 
and the $\widetilde\eta$-invariants of 2-knots 
}
\author{Eiji Ogasa}
\address{
Department of Physics,   University of Tokyo 
Hongo, Tokyo 113, JAPAN 
}
\email{
 ogasa@hep-th.phys.u-tokyo.ac.jp
}

\begin{document}
\begin{abstract} 
We discuss the ribbon-move for 2-knots, which is a local move. 
Let $K$ and $K'$ be 2-knots. Then we have:
Suppose that $K$ and $K'$ are ribbon-move equivalent. 

\noindent
(1) Let ${\mathrm {Tor}} H_1(\widetilde X_K; {\Z})$ 
(resp. ${\mathrm {Tor}} H_1(\widetilde X_{K'}; {\Z})$) 
be the $\Z$-torsion submodule of 
the Alexander module 
$H_1(\widetilde X_K; {\Z})$ (resp. $H_1(\widetilde X_{K'}; {\Z})$).  
Then 
${\mathrm {Tor}} H_1(\widetilde X_K; {\Z})$
is isomorphic to 
${\mathrm {Tor}} H_1(\widetilde X_{K'}; {\Z})$
not only as $\Z$-modules but also as ${\Z}[t,t^{-1}]$-modules. 

\noindent
(2) The Farber-Levine pairing for $K$ is equivalent to that for $K'$. 

\noindent
(3) The set of the values of 
the $\Q/\Z$-valued $\tilde\eta$ invariants for $K$ 
is equivalent to that for $K'$. 
\end{abstract}
\maketitle

\section{Ribbon-moves of 2-knots}\label{Ribbon-moves of 2-knots}

This paper is a sequel of the author's preprint \cite{O2}. 
When the content of this paper and that of \cite{O2} will be published,  
they will be made into one paper.  
For the convenience of the readers, 
this paper includes some parts of \cite{O2}. 
See \cite{O2} for other parts.
The results of this paper include all of the results of \cite{O2} 
and are  stronger than those of \cite{O2}.  
Theorems in Theorem 2.1.(2) is a  new result in this paper. 

In this paper we discuss ribbon-moves of 2-knots.  
In this section we review the definition of ribbon-moves.

An  {\it (oriented) 2-(dimensional) knot}
 is a smooth, oriented submanifold $K$ of $S^4$ 
which is diffeomorphic to the $2$-sphere. 
We say that 2-knots $K_1$ and $K_2$ are {\it equivalent} 
if there exists an orientation preserving diffeomorphism  
$f:$ $S^4$ $\rightarrow$ $S^4$ 
such that $f(K_1)$=$K_2$  and 
that $f | _{K_1}:$ $K_1$ $\rightarrow$ $K_2$ is 
an orientation preserving diffeomorphism.    
Let $id:S^4$ $\rightarrow$ $S^4$ be the identity. 
We say that 2-knots $K_1$ and $K_2$ are {\it identical}  
if  $id(K_1)$=$K_2$ and $id | _{K_1}:K_1\rightarrow K_2$ is 
an orientation preserving diffeomorphism.    

\vskip3mm
\begin{defn}\label{ribbonmove}     
Let $K_1$ and $K_2$ be 2-knots in $S^4$. 
We say that $K_2$ is obtained from $K_1$ by one {\it ribbon-move } 
if there is a 4-ball $B$ embedded in $S^4$ with the following properties.  

\vskip3mm
(1) 
$K_1$ coincides with $K_2$ in $\overline{S^4-B}$. 
This identity map from $\overline{K_1-B}$ to $\overline{K_2-B}$  
is orientation preserving.

(2) 
$B\cap K_1$ is drawn as in Figure 1.1.    
$B\cap K_2$ is drawn as in Figure 1.2.    
\vskip3mm


We regard $B$ as 
(a closed 2-disc)$\times[0,1]\times\{t| -1\leqq t\leqq1\}$.
We put $B_t=$(a closed 2-disc)$\times[0,1]\times\{t \}$.  
Then $B=\cup B_t$. 
In Figure 1.1 
and 1.2, 
we draw $B_{-0.5}, B_{0}, B_{0.5}$ $\subset B$. 
We draw $K_1$ and $K_2$ by the bold line. 
The fine line denotes $\partial B_t$. 
  
$B\cap K_1$ (resp. $B\cap K_2$) is diffeomorphic to 
$D^2\amalg (S^1\times [0,1])$, 
where $\amalg$ denotes the disjoint union.

$B\cap K_1$ has the following properties:  
$B_t\cap K_1$ is empty for $-1\leqq t<0$ and $0.5<t\leqq1$.
$B_0\cap K_1$ is diffeomorphic to 
$D^2\amalg(S^1\times [0,0.3])\amalg(S^1\times [0.7,1])$. 
$B_{0.5}\cap K_1$ is diffeomorphic to $(S^1\times [0.3,0.7])$. 
$B_t\cap K_1$ is diffeomorphic to $S^1\amalg S^1$ for $0<t<0.5$. 
(Here we draw $S^1\times [0,1]$ to have the corner 
in $B_0$ and in $B_{0.5}$. 
Strictly to say, $B\cap K_1$ in $B$ is a smooth embedding 
which is obtained by making the corner smooth naturally.)

$B\cap K_2$ has the following properties:  
$B_t\cap  K_2$ is empty for $-1\leqq t<-0.5$ and $0<t\leqq1$.
$B_0\cap K_2$ is diffeomorphic to 
$D^2\amalg(S^1\times [0, 0.3])\amalg(S^1\times [0.7, 1])$. 
$B_{-0.5}\cap  K_2$ is diffeomorphic to $(S^1\times [0.3, 0.7])$. 
$B_t\cap  K_2$ is diffeomorphic to $S^1\amalg S^1$ for $-0.5<t<0$. 
(Here we draw $S^1\times [0,1]$ to have the corner 
in $B_0$ and in $B_{-0.5}$. 
Strictly to say, $B\cap K_1$ in $B$ is a smooth embedding 
which is obtained by making the corner smooth naturally.)

In Figure 1.1 (resp. 1.2) 
there are an oriented cylinder $S^1\times [0,1]$ 
and an oriented disc $D^2$ as we stated above. 
We do not make any assumption about 
the orientation of the cylinder and the disc. 
The orientation of $B\cap K_1$ (resp. $B\cap K_2$ ) 
coincides with that of the cylinder and that of the disc.

\hskip3cm Figure 1.1.   

\hskip3cm Figure 1.2. 


Suppose that $K_2$ is obtained from $K_1$ by one ribbon-move 
and that $K'_2$ is equivalent to $K_2$.   
Then we also say that $K'_2$ is obtained from $K_1$ 
by one {\it ribbon-move}.   
If $K_1$ is obtained from $K_2$ by one ribbon-move,  
then we also say that $K_2$ is obtained from $K_1$ by one {\it ribbon-move}.   
\end{defn}

\vskip3mm
\begin{defn}\label{ribbonmove2}    
Two 2-knots $K_1$ and $K_2$ are said to be {\it ribbon-move equivalent} 
if there are 2-knots 
$K_1=\bar{K}_1, \bar{K}_2,...,\bar{K}_{r-1},\bar{K}_r=K_2$  
 ($r\in{\N}, p\geq2$) such that 
$\bar{K}_i$ is obtained from $\bar{K}_{i-1}$ $(1< i\leqq r)$ by one ribbon-move. \end{defn}
\vskip3mm

\vskip3mm
\begin{prob}\label{problem}    
Let $K_1$ and $K_2$ be 2-knots.
Find a necessary (resp. sufficient, necessary and sufficient )
condition that $K_1$ and $K_2$ are ribbon-move equivalent.  
\end{prob}
\vskip3mm

In \cite{O} the author proved:  

\vskip3mm
\begin{thm}\label{mu}    
{\rm (\cite{O})}  
(1)If 2-knots $K$ and $K'$ are ribbon-move equivalent, then 
$$\mu(K)=\mu(K'),$$ 

where $\mu(\hskip1.5mm)$ is the $\mu(\hskip1.5mm)$ invarinat of 2-knots 
defined naturally from the $\mu(\hskip1.5mm)$ invariant of 3-manifolds.

\noindent
(2)Let $K_1$ and $K_2$ be 2-knots in $S^4$. 
Suppose that $K_1$ are ribbon-move equivalent to $K_2$.   
Let $W_i$ be an arbitrary Seifert hypersurface for $K_i$. 
Then the torsion part of 
$\{H_1(W_1)\oplus H_1(W_2)\}$ is congruent to 
$G\oplus G$ for a finite abelian group $G$. 

\noindent
(3)Not all 2-knots are ribbon-move equivalent to the trivial 2-knot.   

\noindent
(4)The converse of (1) is not true.  The converse of (2) is not true. 
\end{thm}
\vskip3mm

\section{ Main results }

Let $K$ be a 2-knot. 
Let $\widetilde X^\infty_K$ 
be the canonical infinite cyclic covering space for $K$.  
The homology group $H_1(\widetilde X^\infty_K; {\Z})$ is 
not only a $\Z$-module but also a ${\Z}[t,t^{-1}]$-module.
Let ${\mathrm {Tor}} H_1(\widetilde X^\infty_K; {\Z})$ 
denote the $\Z$-torsion submodule of $H_1(\widetilde X^\infty_K; {\Z})$.  
Then ${\mathrm {Tor}} H_1(\widetilde X^\infty_K; {\Z})$ 
is not only a $\Z$-module but also a ${\Z}[t,t^{-1}]$-module. 

\vskip3mm
\begin{thm}\label{main}    
Let $K$ and $K'$ be 2-knots. 
Suppose that $K$ and $K'$ are ribbon-move equivalent. 
Then there is an isomorphism  
$$c: {\mathrm {Tor}} H_1(\widetilde X^\infty_K; {\Z})\to 
{\mathrm {Tor}} H_1(\widetilde X^\infty_{K'}; {\Z}),$$ 
where the homomorphism $c$ is 
not only one as $\Z$-modules but also one as ${\Z}[t,t^{-1}]$-modules, 
with the following properties. 

\noindent(1) 
Let $x,y\in{\mathrm{Tor}} H_1(\widetilde X^\infty_K;  {\Z})$. 
Then we have 
$${\mathrm {lk}}(x, y)={\mathrm {lk}}(c(x), c(y)),$$ 
where ${\mathrm {lk}}(\quad)$ denotes the Farber-Levine pairing.  
That is, 
the Farber-Levine pairing on 
${\mathrm{Tor}} H_1(\widetilde X^\infty_K;  {\Z})$ is equivalent to 
that on 
${\mathrm{Tor}} H_1(\widetilde X^\infty_{K'};  {\Z})$. 

\noindent(2)
Let $\alpha: H_1(\widetilde{X}^{\infty}_{K}; {\Z})\to{\Z_p}$  
be a homomorphism. 
Note that there is a homomorphism 
\newline
$\alpha': H_1(\widetilde{X}^{\infty}_{K'}; {\Z})\to{\Z_p}$  
such that 
\newline 
$\alpha\vert_{\mathrm Tor}=(\alpha'\vert_{\mathrm Tor})\circ c$. 
Then we have 
$\tilde\eta(K, \alpha)=\tilde\eta(K', \alpha')\in\\Q/\Z$.
That is, 
the set of the values of the $\Q/\Z$-valued $\tilde\eta$ invariants for $K$ 
is equivalent to that for $K'$. 
\end{thm}
\vskip3mm

\noindent{\bf Note.}
Theorem 2.1 (2) is a new result in this paper. 
Theorem 2.1 other than (2) is proved in  \cite{O2}. 

\noindent{\bf Note.}
See [1-17] for the Alexander module, 
the Farber-Levine pairing, 
the $\widetilde\eta$-invariants of 2-knots, 
and the related topics. 


\section{Proof of Theorem \ref{main}(2)}\label{proof}

By using \cite{O2}, we prove Theorem \ref{main}(2).
we complete the proof of  Theorem \ref{main}.

\begin{prop}\label{sagainteger} 
Let $V$ be a closed oriented 3-manifold. 
Let $\alpha, \beta: H_1(V; {\Z})\rightarrow{\Z_p}$ be homomorphisms, 
where $p$ is a natural number greater than one. 
Let $\{t_1,...,t_\mu\}$ be a set of generators of Tor $H_1(V; {\Z})$. 
Suppose that $\alpha(t_i)=\beta(t_i)$ for each $i$. 
Then we have 
$\widetilde\eta(V, \alpha)-\widetilde\eta(V, \beta)\in {\Z}$. 
\end{prop}

\vskip3mm
\noindent{\bf Proof of Proposition \ref{sagainteger}.} 
Recall that the homomorphism $\alpha$ (resp. $\beta$) induces 
a continuous map $V\rightarrow K({\Z_p}, 1)$, 
call it $\alpha$ (resp. $\beta$) again. 
Let $T_i$ be a tubular neighborhood of $t_i$. 
Let $N=T_1\amalg...\amalg T_\mu$. 
Then we can suppose that $\alpha\vert_N=\beta\vert_N$.

Take $V\x[0,3]$. 
Take $N\x[1,2]\subset V\x[0,3]$. 
Define a map 

\noindent
$f:V\x[0,1]\cup N\x[1,2]\cup V\x[2,3]\rightarrow K({\Z_p}, 1)$ 
with the following properties. 

\vskip3mm
\noindent\hskip5mm(1)  
$f\vert_{V\x[0,1]}:V\x[0,1]\rightarrow K({\Z_p}, 1)$ 
is defined by 

\noindent\hskip5mm
$(x, t)\mapsto\alpha(x)$. 

\noindent\hskip5mm(2)
$f\vert_{V\x[2,3]}:V\x[2,3]\rightarrow K({\Z_p}, 1)$ 
is defined by 

\noindent\hskip5mm
$(x, t)\mapsto\beta(x)$.

\noindent\hskip5mm(3)
$f\vert_{N\x[1,2]}:N\x[1,2]\rightarrow K({\Z_p}, 1)$ 
is defined by 

\noindent\hskip5mm
$(x, t)\mapsto\alpha(x)(=\beta(x))$. 

\noindent
Let $X=V\x[0,1]\cup N\x[1,2]\cup V\x[2,3]$.  
Let $M=\partial X-(V\x\{0\})-(V\x\{3\})$. 
Let $\gamma=f\vert_M$. 
Then we have 

\vskip3mm
\hskip7mm
$\partial(V\x[0,1]\cup N\x[1,2]\cup V\x[2,3], f)$
  
\noindent\hskip7mm  
$=(V\x\{0\}, f\vert_{V\x\{0\}})\amalg 
(V\x\{3\}, f\vert_{V\x\{3\}})\amalg (M, f\vert_M)$ 

\noindent\hskip7mm  
$=(V, \alpha)\amalg -(V, \beta)\amalg (M, \gamma)$.   
\vskip3mm


\noindent
Hence we have  
$\widetilde\eta(V, \alpha)-\widetilde\eta(V, \beta)=\widetilde\eta(M, \gamma)$. 
From now we investigate $M$. 
Note that 
$M=\partial(\overline{(V-N)\times[1,2]})$.

\vskip3mm 
\begin{cla}\label{VN}
We have Tor$H_i(V-N; {\Z})\cong0$ for each $i$. 
\end{cla} 
\vskip3mm 

\noindent{\bf Proof of Claim \ref{VN}.} 
Since $V-N$ is a compact oriented 3-manifold with boundary 
and $\partial(V-N)\neq\phi$, 
there is a handle decomposition of one 0-handle, 1-handles, and 2-handles. 
By using this handle decomposition, 
we can calculate $H_i(V-N; {\Z})$. 
Hence we have Tor$H_i(V-N; {\Z})\cong0$ for $i\neq1$. 

By the Poincar\'e duality, the universal coefficient theorem, 
and the excision, we have 
Tor$H_1(V-N; {\Z})\cong$ 
Tor$H_1(V-N, \partial(V-N); {\Z})\cong$ 
Tor$H_1(V, N; {\Z})$.  
Consider the Mayer-Vietoris exact sequence: 

$$
H_i(N; {\Z})\stackrel{\iota}\rightarrow 
H_i(V; {\Z})\stackrel{q}\rightarrow 
H_i(V, N; {\Z}). 
$$

\noindent
By the definition of $N$, it holds that 

\[
H_i(N; {\Z})\cong
\left\{
\begin{array}{ll}
0&\mbox{if $i\neq0,1$}\\
{\Z^\mu}&\mbox{if $i=0,1$}
\end{array}
\right.
\]

\noindent 
Consider 
$H_1(N; {\Z})\stackrel{\iota}\rightarrow H_1(V; {\Z})$. 
By the definition of the inclusion $N\rightarrow V$, 
it holds that Im$\iota\cong$ Tor$H_1(V; {\Z})$.
Furthermore, we have 
Im $\{H_1(V, N; {\Z})\stackrel{\partial}\rightarrow H_0(N; {\Z})\}$ 
is torsion free. 
Hence $H_1(V, N; {\Z})$ is torsion free. 
Hence $H_1(V- N; {\Z})$ is torsion free. 
This completes the proof of Claim \ref{VN}.

\vskip3mm 
\begin{cla}\label{MM}
We have Tor$H_1(M; {\Z})\cong0$. 
\end{cla} 

\vskip3mm 
\noindent{\bf Proof of Claim \ref{MM}.} 
Recall that
$M=\partial(\overline{(V-N)\times[1,2]})$. 
Consider the Mayer-Vietoris exact sequence:  

$$
H_i(M; {\Z})
\stackrel{\iota}\rightarrow 
H_i((V-N)\x[1,2]; {\Z})
\stackrel{q}\rightarrow 
H_i((V-N)\x[1,2], M; {\Z}) 
$$

\noindent
and consider the following part:  

$$
H_2((V-N)\x[1,2]; {\Z})
\stackrel{q}\rightarrow 
H_2((V-N)\x[1,2], M; {\Z}) 
\stackrel{\partial}\rightarrow 
H_1(M; {\Z})
\stackrel{\iota}\rightarrow 
H_1((V-N)\x[1,2]; {\Z}). 
$$

\noindent
By applying the Poincar\'e duality 
and the universal coefficient theorem, we can obtain 
\newline 
$H_i((V-N)\x[1,2]; {\Z})$ 
and 
$H_i((V-N)\x[1,2], M; {\Z})$ 
\newline
by using $H_i(V-N; {\Z})$. 
 By Claim \ref{VN}, 
 $H_1(V- N; {\Z})$ is torsion free. 
\newline 
Hence we have 
(1)$H_2((V-N)\x[1,2], M; {\Z})$ 
is torsion free, put it to be 
${\Z^\beta}$, and 
(2) $H_1((V-N)\x[1,2]; {\Z})$) 
is torsion free, put it to be ${\Z^\gamma}$. 

\noindent
Let $x$ be any 2-cycle in $(V-N)\x[1,2]$. 
Then we can suppose that 
$x\subset (V-N)\x\{1\}\subset M$. 
Hence  
$H_2((V-N)\x[1,2]; {\Z})\stackrel{q}\rightarrow 
H_2((V-N)\x[1,2], M; {\Z})$ is the zero map.  
Hence we have an exact sequence: 

$$
0
\rightarrow 
{\Z^\beta}
\rightarrow 
H_1(M; {\Z})
\rightarrow 
{\Z^\gamma}.
$$

\noindent
Hence $H_1(M; {\Z})$ is torsion free. 
This completes the proof of Claim \ref{MM}. 
\vskip3mm 

\begin{cla}\label{integer} 
Let $M$ be a closed oriented 3-manifold. 
Suppose that $H_1(M; {\Z})$ is torsion free. 
Let $\gamma$ be a homomorphism $H_1(M; {\Z})\rightarrow{\Z_p}$. 
Then we have $\widetilde\eta(M, \gamma)\in {\Z}$. 
\end{cla}

\vskip3mm\noindent{\bf Proof of Claim \ref{integer}.} 
By an elementary discussion on homomorphisms, we have 
the following commutative diagram. 

$$
\begin{array}{ccc}
{\Z}&\stackrel{\gamma'}\rightarrow&{\Z_p}\\ 
\uparrow_{\zeta} &\nearrow_\gamma &\\
H_1(M;{\Z})&&
\end{array}, 
$$

\noindent
where the above homomorphism $H_1(M;{\Z})\to{\Z}$ is called $\zeta$.
Then $\gamma'\circ\zeta=\gamma$. Then we can regard  

\vskip3mm
$\zeta\in$Hom$(H_1(M;{\Z}),{\Z})$

$\hskip1cm\cong 
H^1(M;{\Z})$

$\hskip1cm\cong  
\{\mathrm{homotopy\hskip1mm classes\hskip1mm of\hskip1mm maps}\hskip1mm  
M\to K({\Z}, 1) \}$ 

\vskip3mm
\noindent and 
\vskip3mm

$\gamma'\circ\zeta\in$Hom$(H_1(M;{\Z}),{\Z_p})$

$\hskip1cm\cong H^1(M;{\Z_p})$

$\hskip1cm\cong
\{\mathrm{homotopy\hskip1mm classes\hskip1mm of\hskip1mm maps}
\hskip1mm    
M\to K({\Z_p},1)\}.$ 
\vskip3mm

\noindent 
The above commutative diagram induces 
the following commutative diagram. 

$$
\begin{array}{ccc}
K({\Z},1)&\stackrel{\gamma'}\rightarrow&K({\Z_p},1)\\ 
\uparrow_\zeta &\nearrow_\gamma &\\
M&&
\end{array}
$$

\noindent
The above commutative diagram induces the following homomorphism. 

$$
\begin{array}{cccc}
&\Omega^3(K({\Z},1))&\stackrel{\gamma'}\rightarrow&\Omega^{3}(K({\Z_p},1))\\
&\rotatebox[origin=c]{90}{$\in $}& &\rotatebox[origin=c]{90}{$\in $}\\
&[(M, \zeta)] & \mapsto& [(M, \gamma'\circ\zeta)]  
\end{array}
$$

\noindent
By bordism theory, we have  
$\Omega^3(K({\Z}, 1))\cong\Omega^3(S^1)\cong0$. 
Therefore  we have 
$[(M,\gamma'\circ\zeta)]=0$. 
Therefore we have 
$(M, \gamma'\circ\zeta)=\partial(W,\tau)$ 
for a compact oriented 4-manifold $W$ and a continuous map $\tau$.   
Hence we have 
$\widetilde\eta(M,\gamma'\circ\zeta )=
\frac{1}{1}(\overline\sigma(W, \tau)-\sigma(W))$. 
Hence we have $\widetilde\eta(M,\gamma'\circ\zeta)\in{\Z}$. 
Hence we have $\widetilde\eta(M,\gamma)\in{\Z}$. 
This completes the proof of Claim \ref{integer}. 
\vskip3mm 

Proposition \ref{sagainteger} holds by Claim \ref{MM}, \ref{integer}.  
\vskip3mm

\noindent{\bf Note.}
Furthermore we could prove:     
Let $X$ and $Y$ be closed oriented 3-manifolds. 
Suppose that there is an isomorphism 
$c: {\mathrm {Tor}} H_1(X;{\Z})\to 
{\mathrm {Tor}} H_1(Y;{\Z})$ 
 such that 
${\mathrm {lk}}(\epsilon,\delta)=
{\mathrm {lk}}(c(\epsilon),c(\delta))$
for any pair of elements  $\epsilon,\delta\in{\mathrm {Tor}} H_1(X;{\Z})$. 
Let $f:X\to B\Z_p$ and $g:Y\to B\Z_p$ be continuous maps 
such that the homomorphism 
$f_\star\vert_\mathrm{Tor}:\mathrm{Tor}H_1(X;\Z)\to\Z_p$ 
and 
$g_\star\vert_\mathrm{Tor}:\mathrm{Tor}H_1(Y;\Z)\to\Z_p$ 
satisfy 
$f_\star=c\circ g_\star$.
Then we have 
$\widetilde\eta(X, f)=\widetilde\eta(Y,g)$ mod$\Z$.


\vskip3mm
Recall we prove the following in \cite{O2}\cite{O1}.

\vskip3mm 
\begin{prop}\label{torsiondiagram} 
Suppose that $K_<$ is obtained from $K_>$ by one ribbon move. 
Then we have the following. 

\noindent 
(1) 
There are Seifert hypersurfaces $V_<$ for $K_<$ and 
$V_>$ for $K_>$ such that $V_<$ is diffeomorphic to $V_>$. 

\noindent
(2) There are a lift $V_{<,\xi}$ of $V_<$ in $\widetilde{X}^{\infty}_{K_<}$ 
    a lift $V_{>,\xi}$ of $V_>$ in $\widetilde{X}^{\infty}_{K_>}$, 
    and a diffeomorphism $h: V_{<,\xi}\rightarrow V_{>,\xi}$ such that 
    we have the following commutative diagram.       

\[\begin{array}{ccc}
{\mathrm {Tor}} H_1(V_{<,\xi}; {\Z}) & 
\stackrel{h_*\vert_{\mathrm Tor}, \cong}\to  &
{\mathrm {Tor}} H_1(V_{>,\xi}; {\Z}) 
\\ 
\downarrow ^{\iota_<\vert_{\mathrm {Tor}}, {\mathrm{onto}}}&
&\downarrow^{\iota_>\vert_{\mathrm {Tor}}, {\mathrm{onto}}} 
\\ 
{\mathrm {Tor}} H_1(\widetilde{X}^{\infty}_{K_<}; {\Z}) & 
\stackrel{c,\cong}\to &
{\mathrm {Tor}} H_1(\widetilde{X}^{\infty}_{K_>}; {\Z}), 
\\ 
\end{array}\] 

\noindent
(3) 
Let $x,y\in{\mathrm Tor} H_1(\widetilde{X}^{\infty}_{K_<}; {\Z})$.  
Let ${\mathrm {lk}}(\quad, \quad)$ denote the Farber-Levine pairing.  
Then ${\mathrm {lk}}(x, y)={\mathrm {lk}}(c(x), c(y))$.  
\end{prop}

We prove:

\vskip3mm
\begin{prop}\label{uptildeeta} 
We can add the following condition (4) to Proposition \ref{torsiondiagram}. 

\noindent
(4)
Let $\alpha_<: H_1(\widetilde{X}^{\infty}_{K_<}; {\Z})\to{\Z_p}$  
be a homomorphism. 
Note that there is a homomorphism 
\newline
$\alpha_>: H_1(\widetilde{X}^{\infty}_{K_>}; {\Z})\to{\Z_p}$  
such that 
\newline 
$\alpha_<\vert_{\mathrm Tor}=(\alpha_>\vert_{\mathrm Tor})\circ c$. 
Then we have 
$\tilde\eta(K_<, \alpha_<)=\tilde\eta(K_>, \alpha_>)$. 
That is, 
the set of the values of the $\widetilde\eta$-invariants of 
the 2-knot $K_<$ 
coincides with that of $K_>$ as above.   
\end{prop}

\vskip3mm
\noindent{\bf Proof of Proposition \ref{uptildeeta}.} 
By the definition of the $\tilde\eta$-invariants of 2-knots, we have 
\newline 
$\tilde\eta(K_<, \alpha_<)=\tilde\eta(V_<, \alpha_<\circ\iota_<)$ 
mod ${\Z}$ and 
\newline 
$\tilde\eta(K_>, \alpha_>)=\tilde\eta(V_>, \alpha_>\circ\iota_>)$ 
mod ${\Z}$. 
By using the commutative diagram in Proposition \ref{torsiondiagram}, 
we have  
\newline 
$(\alpha_<\vert{\mathrm Tor})\circ(\iota_<\vert{\mathrm Tor})
=(\alpha_>\vert{\mathrm Tor})\circ c\circ(\iota_<\vert{\mathrm Tor})
=(\alpha_>\vert{\mathrm Tor})\circ(\iota_>\vert{\mathrm Tor})\circ h
: {\mathrm Tor}H_1(V_<; {\Z})\to{\Z_p}$. 
\newline
Hence we have 

\noindent$\tilde\eta(V_<, \alpha_<\circ \iota_<)$

\noindent$=\tilde\eta(V_<, \alpha_>\circ \iota_>\circ h)$

\noindent$=\tilde\eta(h(V_<), \alpha_>\circ \iota_>)$

\noindent$=\tilde\eta(V_>, \alpha_>\circ \iota_>)$. 

This completes the proof of Proposition \ref{uptildeeta}.


\vskip3mm
By using Proposition \ref{uptildeeta} some times, 
Theorem \ref{main} holds.


\noindent{\bf Acknowledgment.}
The author would like to thank Prof. Kent Orr 
for correcting his English and for the valuable discussion.

\np

\pagestyle{empty}

\unitlength 0.1in
\begin{picture}(56.10,47.00)(8.50,-47.60)
%
\special{pn 8}%
\special{ar 3510 320 560 250  0.0000000 6.2831853}%
%
\special{pn 20}%
\special{pa 3660 3950}%
\special{pa 3660 2730}%
\special{fp}%
%
\special{pn 20}%
\special{pa 3340 3970}%
\special{pa 3340 2740}%
\special{fp}%
%
\special{pn 20}%
\special{pa 3350 4000}%
\special{pa 3360 3970}%
\special{pa 3384 3949}%
\special{pa 3413 3935}%
\special{pa 3444 3926}%
\special{pa 3475 3921}%
\special{pa 3507 3920}%
\special{pa 3539 3923}%
\special{pa 3570 3930}%
\special{pa 3601 3940}%
\special{pa 3627 3958}%
\special{pa 3647 3983}%
\special{pa 3648 4014}%
\special{pa 3630 4040}%
\special{pa 3603 4058}%
\special{pa 3573 4070}%
\special{pa 3542 4077}%
\special{pa 3510 4080}%
\special{pa 3478 4079}%
\special{pa 3447 4075}%
\special{pa 3416 4066}%
\special{pa 3387 4053}%
\special{pa 3362 4032}%
\special{pa 3350 4003}%
\special{pa 3350 4000}%
\special{sp}%
%
\special{pn 20}%
\special{pa 3340 360}%
\special{pa 3340 1580}%
\special{fp}%
%
\special{pn 20}%
\special{pa 3660 340}%
\special{pa 3660 1570}%
\special{fp}%
%
\special{pn 20}%
\special{ar 3500 310 150 80  0.0000000 6.2831853}%
%
\special{pn 8}%
\special{ar 3510 4080 560 250  0.0000000 6.2831853}%
%
\special{pn 8}%
\special{pa 4080 330}%
\special{pa 4080 4080}%
\special{fp}%
%
\special{pn 8}%
\special{pa 2950 340}%
\special{pa 2950 4080}%
\special{fp}%
%
\special{pn 8}%
\special{ar 1410 330 560 250  0.0000000 6.2831853}%
%
\special{pn 8}%
\special{ar 1410 4090 560 250  0.0000000 6.2831853}%
%
\special{pn 8}%
\special{pa 1980 340}%
\special{pa 1980 4090}%
\special{fp}%
%
\special{pn 8}%
\special{pa 850 350}%
\special{pa 850 4090}%
\special{fp}%
%
\special{pn 8}%
\special{ar 5890 310 560 250  0.0000000 6.2831853}%
%
\special{pn 8}%
\special{ar 5890 4070 560 250  0.0000000 6.2831853}%
%
\special{pn 8}%
\special{pa 6460 320}%
\special{pa 6460 4070}%
\special{fp}%
%
\special{pn 8}%
\special{pa 5330 330}%
\special{pa 5330 4070}%
\special{fp}%
%
\special{pn 20}%
\special{ar 3500 1590 150 80  0.0000000 6.2831853}%
%
\special{pn 20}%
\special{ar 3500 2710 150 80  0.0000000 6.2831853}%
%
\special{pn 20}%
\special{ar 3520 2070 560 250  0.0000000 6.2831853}%
%
\special{pn 8}%
\special{pa 3660 1600}%
\special{pa 5730 1600}%
\special{dt 0.045}%
\special{pa 5730 1600}%
\special{pa 5729 1600}%
\special{dt 0.045}%
%
\special{pn 8}%
\special{pa 3680 2720}%
\special{pa 5750 2720}%
\special{dt 0.045}%
\special{pa 5750 2720}%
\special{pa 5749 2720}%
\special{dt 0.045}%
%
\special{pn 20}%
\special{ar 5850 1600 150 80  0.0000000 6.2831853}%
%
\special{pn 20}%
\special{ar 5850 2730 150 80  0.0000000 6.2831853}%
%
\special{pn 20}%
\special{pa 5690 1610}%
\special{pa 5690 2700}%
\special{fp}%
%
\special{pn 20}%
\special{pa 6020 1640}%
\special{pa 6020 2730}%
\special{fp}%
\put(12.0000,-46.0000){\makebox(0,0)[lb]{t=-0.5}}%
\put(32.0000,-46.0000){\makebox(0,0)[lb]{t=0}}%
\put(56.0000,-46.0000){\makebox(0,0)[lb]{t=0.5}}%
\put(30.3000,-49.3000){\makebox(0,0)[lb]{Figure 1.1}}%
%
\special{pn 8}%
\special{ar 3510 320 560 250  0.0000000 6.2831853}%
%
\special{pn 8}%
\special{pa 2950 340}%
\special{pa 2950 4080}%
\special{fp}%
%
\special{pn 8}%
\special{pa 4080 330}%
\special{pa 4080 4080}%
\special{fp}%
%
\special{pn 8}%
\special{ar 3510 4080 560 250  0.0000000 6.2831853}%
%
\special{pn 20}%
\special{ar 3500 310 150 80  0.0000000 6.2831853}%
%
\special{pn 20}%
\special{pa 3660 340}%
\special{pa 3660 1570}%
\special{fp}%
%
\special{pn 20}%
\special{pa 3340 360}%
\special{pa 3340 1580}%
\special{fp}%
%
\special{pn 20}%
\special{pa 3350 4000}%
\special{pa 3360 3970}%
\special{pa 3384 3949}%
\special{pa 3413 3935}%
\special{pa 3444 3926}%
\special{pa 3475 3921}%
\special{pa 3507 3920}%
\special{pa 3539 3923}%
\special{pa 3570 3930}%
\special{pa 3601 3940}%
\special{pa 3627 3958}%
\special{pa 3647 3983}%
\special{pa 3648 4014}%
\special{pa 3630 4040}%
\special{pa 3603 4058}%
\special{pa 3573 4070}%
\special{pa 3542 4077}%
\special{pa 3510 4080}%
\special{pa 3478 4079}%
\special{pa 3447 4075}%
\special{pa 3416 4066}%
\special{pa 3387 4053}%
\special{pa 3362 4032}%
\special{pa 3350 4003}%
\special{pa 3350 4000}%
\special{sp}%
%
\special{pn 20}%
\special{pa 3340 3970}%
\special{pa 3340 2740}%
\special{fp}%
%
\special{pn 20}%
\special{pa 3660 3950}%
\special{pa 3660 2730}%
\special{fp}%
%
\special{pn 8}%
\special{ar 1410 330 560 250  0.0000000 6.2831853}%
%
\special{pn 8}%
\special{ar 1410 4090 560 250  0.0000000 6.2831853}%
%
\special{pn 8}%
\special{pa 1980 340}%
\special{pa 1980 4090}%
\special{fp}%
%
\special{pn 8}%
\special{pa 850 350}%
\special{pa 850 4090}%
\special{fp}%
%
\special{pn 8}%
\special{ar 5890 310 560 250  0.0000000 6.2831853}%
%
\special{pn 8}%
\special{ar 5890 4070 560 250  0.0000000 6.2831853}%
%
\special{pn 8}%
\special{pa 6460 320}%
\special{pa 6460 4070}%
\special{fp}%
%
\special{pn 8}%
\special{pa 5330 330}%
\special{pa 5330 4070}%
\special{fp}%
\end{picture}%

\np
\unitlength 0.1in
\begin{picture}(56.10,47.10)(8.50,-47.70)
%
\special{pn 8}%
\special{ar 3510 320 560 250  0.0000000 6.2831853}%
%
\special{pn 20}%
\special{pa 3660 3950}%
\special{pa 3660 2730}%
\special{fp}%
%
\special{pn 20}%
\special{pa 3340 3970}%
\special{pa 3340 2740}%
\special{fp}%
%
\special{pn 20}%
\special{pa 3350 4000}%
\special{pa 3360 3970}%
\special{pa 3384 3949}%
\special{pa 3413 3935}%
\special{pa 3444 3926}%
\special{pa 3475 3921}%
\special{pa 3507 3920}%
\special{pa 3539 3923}%
\special{pa 3570 3930}%
\special{pa 3601 3940}%
\special{pa 3627 3958}%
\special{pa 3647 3983}%
\special{pa 3648 4014}%
\special{pa 3630 4040}%
\special{pa 3603 4058}%
\special{pa 3573 4070}%
\special{pa 3542 4077}%
\special{pa 3510 4080}%
\special{pa 3478 4079}%
\special{pa 3447 4075}%
\special{pa 3416 4066}%
\special{pa 3387 4053}%
\special{pa 3362 4032}%
\special{pa 3350 4003}%
\special{pa 3350 4000}%
\special{sp}%
%
\special{pn 20}%
\special{pa 3340 360}%
\special{pa 3340 1580}%
\special{fp}%
%
\special{pn 20}%
\special{pa 3660 340}%
\special{pa 3660 1570}%
\special{fp}%
%
\special{pn 20}%
\special{ar 3500 310 150 80  0.0000000 6.2831853}%
%
\special{pn 8}%
\special{ar 3510 4080 560 250  0.0000000 6.2831853}%
%
\special{pn 8}%
\special{pa 4080 330}%
\special{pa 4080 4080}%
\special{fp}%
%
\special{pn 8}%
\special{pa 2950 340}%
\special{pa 2950 4080}%
\special{fp}%
%
\special{pn 8}%
\special{ar 1410 330 560 250  0.0000000 6.2831853}%
%
\special{pn 8}%
\special{ar 1410 4090 560 250  0.0000000 6.2831853}%
%
\special{pn 8}%
\special{pa 1980 340}%
\special{pa 1980 4090}%
\special{fp}%
%
\special{pn 8}%
\special{pa 850 350}%
\special{pa 850 4090}%
\special{fp}%
%
\special{pn 8}%
\special{ar 5890 310 560 250  0.0000000 6.2831853}%
%
\special{pn 8}%
\special{ar 5890 4070 560 250  0.0000000 6.2831853}%
%
\special{pn 8}%
\special{pa 6460 320}%
\special{pa 6460 4070}%
\special{fp}%
%
\special{pn 8}%
\special{pa 5330 330}%
\special{pa 5330 4070}%
\special{fp}%
%
\special{pn 20}%
\special{ar 3500 1590 150 80  0.0000000 6.2831853}%
%
\special{pn 20}%
\special{ar 3500 2710 150 80  0.0000000 6.2831853}%
%
\special{pn 20}%
\special{ar 3520 2070 560 250  0.0000000 6.2831853}%
\put(12.0000,-46.0000){\makebox(0,0)[lb]{t=-0.5}}%
\put(32.0000,-46.0000){\makebox(0,0)[lb]{t=0}}%
\put(56.0000,-46.0000){\makebox(0,0)[lb]{t=0.5}}%
%
\special{pn 8}%
\special{ar 3510 320 560 250  0.0000000 6.2831853}%
%
\special{pn 8}%
\special{pa 2950 340}%
\special{pa 2950 4080}%
\special{fp}%
%
\special{pn 8}%
\special{pa 4080 330}%
\special{pa 4080 4080}%
\special{fp}%
%
\special{pn 8}%
\special{ar 3510 4080 560 250  0.0000000 6.2831853}%
%
\special{pn 20}%
\special{ar 3500 310 150 80  0.0000000 6.2831853}%
%
\special{pn 20}%
\special{pa 3660 340}%
\special{pa 3660 1570}%
\special{fp}%
%
\special{pn 20}%
\special{pa 3340 360}%
\special{pa 3340 1580}%
\special{fp}%
%
\special{pn 20}%
\special{pa 3350 4000}%
\special{pa 3360 3970}%
\special{pa 3384 3949}%
\special{pa 3413 3935}%
\special{pa 3444 3926}%
\special{pa 3475 3921}%
\special{pa 3507 3920}%
\special{pa 3539 3923}%
\special{pa 3570 3930}%
\special{pa 3601 3940}%
\special{pa 3627 3958}%
\special{pa 3647 3983}%
\special{pa 3648 4014}%
\special{pa 3630 4040}%
\special{pa 3603 4058}%
\special{pa 3573 4070}%
\special{pa 3542 4077}%
\special{pa 3510 4080}%
\special{pa 3478 4079}%
\special{pa 3447 4075}%
\special{pa 3416 4066}%
\special{pa 3387 4053}%
\special{pa 3362 4032}%
\special{pa 3350 4003}%
\special{pa 3350 4000}%
\special{sp}%
%
\special{pn 20}%
\special{pa 3340 3970}%
\special{pa 3340 2740}%
\special{fp}%
%
\special{pn 20}%
\special{pa 3660 3950}%
\special{pa 3660 2730}%
\special{fp}%
%
\special{pn 8}%
\special{ar 1410 330 560 250  0.0000000 6.2831853}%
%
\special{pn 8}%
\special{ar 1410 4090 560 250  0.0000000 6.2831853}%
%
\special{pn 8}%
\special{pa 1980 340}%
\special{pa 1980 4090}%
\special{fp}%
%
\special{pn 8}%
\special{pa 850 350}%
\special{pa 850 4090}%
\special{fp}%
%
\special{pn 8}%
\special{ar 5890 310 560 250  0.0000000 6.2831853}%
%
\special{pn 8}%
\special{ar 5890 4070 560 250  0.0000000 6.2831853}%
%
\special{pn 8}%
\special{pa 6460 320}%
\special{pa 6460 4070}%
\special{fp}%
%
\special{pn 8}%
\special{pa 5330 330}%
\special{pa 5330 4070}%
\special{fp}%
%
\special{pn 8}%
\special{pa 3320 1590}%
\special{pa 1590 1590}%
\special{dt 0.045}%
\special{pa 1590 1590}%
\special{pa 1591 1590}%
\special{dt 0.045}%
%
\special{pn 8}%
\special{pa 3330 2740}%
\special{pa 1620 2740}%
\special{dt 0.045}%
\special{pa 1620 2740}%
\special{pa 1621 2740}%
\special{dt 0.045}%
%
\special{pn 20}%
\special{ar 1420 1590 150 80  0.0000000 6.2831853}%
%
\special{pn 20}%
\special{ar 1420 2720 150 80  0.0000000 6.2831853}%
%
\special{pn 20}%
\special{pa 1260 1600}%
\special{pa 1260 2690}%
\special{fp}%
%
\special{pn 20}%
\special{pa 1590 1630}%
\special{pa 1590 2720}%
\special{fp}%
\put(30.0000,-49.4000){\makebox(0,0)[lb]{Figure 1.2}}%
\end{picture}%


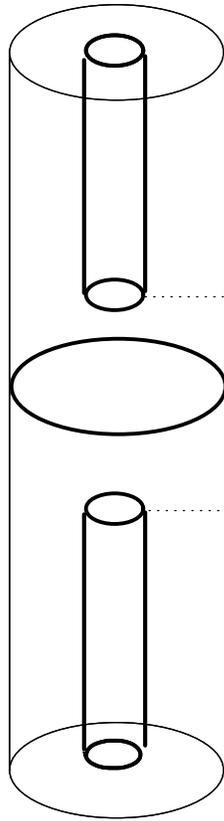
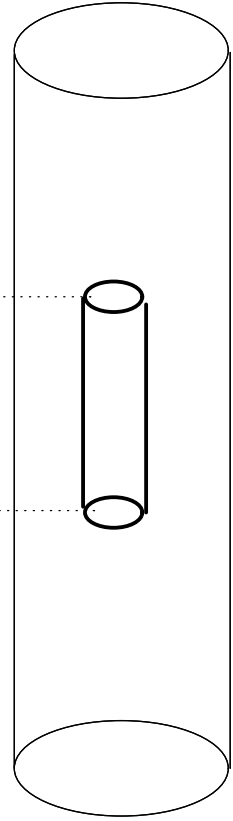
\end{document}